\documentclass[11pt]{amsart}
\usepackage[a4paper]{geometry}
\usepackage{color}
\usepackage{longtable}

\usepackage{lipsum}
\newcommand{\mylongtitle}[1]{%
  \ifodd\value{page}%
    \protect\parbox{0.9\linewidth}{#1}\hfill%
  \else%
    \hfill\protect\parbox{0.9\linewidth}{#1}%
  \fi%
}

 \newtheorem{thm}{Theorem}[section]
 
 \newtheorem{lem}[thm]{Lemma}
 
 \theoremstyle{definition}
 
 \theoremstyle{remark}
 
 \numberwithin{equation}{section}

\begin{document}
\title {\textbf{ \mylongtitle{On Biconservative Lorentz Hypersurface with non-diagonalizable shape operator}}}
\author{ Deepika}

\begin{abstract}
In this paper, we obtain some properties of biconservative Lorentz
hypersurface $M_{1}^{n}$ in $E_{1}^{n+1}$ having shape operator with
complex eigenvalues. We prove that every biconservative Lorentz
hypersurface $M_{1}^{n}$ in $E_{1}^{n+1}$ whose shape operator has
complex eigenvalues with at most five distinct principal curvatures
has constant mean curvature. Also, we investigate such type of
hypersurface with constant length of second fundamental form having
six distinct principal curvatures.
\\
\\
\textbf{AMS2000 MSC Codes:} 53D12, 53C40, 53C42
\end{abstract}

\maketitle
\section{\textbf{Introduction}}

A hypersurface $M^{n}$ in a Riemannian manifold $N^{n+1}$ is called
{\it biconservative} if

\begin{equation}\label{1.1}
 2S (\rm grad \emph{H})+ n \emph{H} \hspace{.1 cm} grad \emph{H} = 2 \emph{H}\hspace{.1 cm}
 Ricci^{N}(\xi)^{\top},
\end{equation}
where $S$ is the shape operator, $H$ is the mean curvature function
and $Ricci^{N}(\xi)^{\top}$ is the tangent component of the Ricci
curvature of $N$ in the direction of the unit normal $\xi$ of
$M^{n}$ in $N^{n+1}$.

In this paper we consider the biconservative Lorentzian hypersurface
$M^{n}_{1}$ in pseudo-Euclidean space $E^{n+1}_{1}$. In this case
(\ref{1.1}) becomes

\begin{equation}\label{1.2}
 2S (\rm grad \emph{H})+ n \emph{H} \hspace{.1 cm} grad \emph{H} = 0,
\end{equation}

From (\ref{1.2}), it is obvious that hypersurfaces with constant
mean curvature are always biconservative. Now, there arise a
question whether there exist biconservative hypersurfaces which are
not of constant mean curvature, known as proper biconservative
\cite{bic10}. In \cite{bic9} and \cite{bic11}, proper biconservative
surfaces in $E^3$ have been classified by proving that they must be
surfaces of revolution. Therefore, it will be interesting to study
the existence/non-existence of proper biconservative hypersurfaces
in pseudo-Euclidean space.

The concept of biconservative hypersurfaces have been studied by
several geometers in [3, 4, 8-12]. The first result on
biconservative hypersurfaces was obtained by T. Hasanis and T.
Vlachos in \cite{bic11}, where biconservative hypersurfaces are
called as H-hypersurfaces. In \cite{bic9}, R. Caddeo et al.
introduced the notion of biconservative and proved that a
biconservative surface in Euclidean 3-space is  either a surface of
constant mean curvature or a surface of revolution (cf.
\cite{bic11}, \cite{bic12}). In \cite{bic3}, the authors proved that
a $\delta(2)-$ideal biconservative hypersurface in Euclidean space
$E^{n}$ $(n \geq 3)$ is either minimal or a spherical hypercylinder.
In \cite{bic10}, Montaldo et al. studied proper $SO(p+1)\times
SO(p+1)$-invariant biconservative hypersurfaces and proper
$SO(p+1)$-invariant biconservative hypersurfaces in Euclidean space
$E^{n}$. Also, Fectu et al. classified biconservative surfaces in
$S^{n}\times R$ and $H^{n}\times R$ in \cite{bic4}. Recently, in
\cite{bic8}, Turgay obtained complete classification of
H-hypersurfaces with three distinct principal curvatures in
Euclidean spaces.

 Our goal is to investigate the nature of
biconservative Lorentz hypersurface $M_{1}^{n}$ in $E_{1}^{n+1}$
whose shape operator has complex eigenvalues. In Section 3, we
obtain some properties of biconservative Lorentz hypersurface in
$E_{1}^{n+1}$ having complex eigenvalues. In Section 4, we study
such type of biconservative Lorentz hypersurface in $E_{1}^{n+1}$
with at most five distinct principal curvatures and concluded the
following result:
\begin{thm}\label{Theorem1.1}
Let $M_{1}^{n}$ in $E^{n+1}_{1}$ be a biconservative Lorentz
hypersurface having non-diagonalizable shape operator with complex
eigenvalues and with at most five distinct principal curvatures.
Then $M_{1}^{n}$ has constant mean curvature.
\end{thm}

In Section 5, we investigate biconservative hypersurface with
constant length of second fundamental form and with six distinct
eigenvalues and obtained the following result:

\begin{thm}\label{Theorem1.2}
Let $M_{1}^{n}$ in $E^{n+1}_{1}$ be a biconservative Lorentz
hypersurface with constant length of second fundamental form and
whose shape operator has complex eigenvalues with six distinct
principal curvatures. Then $M_{1}^{n}$ has constant mean curvature.
\end{thm}

The study of biconservative hypersurfaces is also relevant for the
study of biharmonic hypersurfaces satisfying $\triangle \vec {H} =
0$, where $\triangle$ is a Laplacian operator. It can be seen that
equation (\ref{1.2}) is the tangential component of $\triangle \vec
{H} = 0$, so biharmonic hypersurfaces are always biconservative [2,
12]. Thus, biconservative hypersurfaces form a much larger family of
hypersurfaces including biharmonic hypersurfaces. Recently, the
author has proved that every Lorentz hypersurface in $E^{n+1}_{1}$
satisfying $\triangle \vec {H}= \alpha \vec {H}$ and having complex
eigenvalues with at most four distinct principal curvatures has
constant mean curvature \cite{bic5}. So, it will be interesting to
investigate the biconservative Lorentzian hypersurfaces having
complex eigenvalues as a natural generalization and extension of the
results obtained in \cite{bic5}.

\vspace{.5 cm}
\section{\textbf{Preliminaries}}

Let ($M_{1}^{n}, g$) be an $n$-dimensional  biconservative Lorentz
hypersurface of an $n+1$-dimensional pseudo-Euclidean space
($E^{n+1}_{1}, \overline g$) and $g = \overline g_{|M_{1}^{n}}$. We
denote by $\xi$ the unit normal vector to $M_{1}^{n}$ with
$\overline g(\xi, \xi)= 1$.

  Let $\overline\nabla $ and $ \nabla $ denote the linear connections of $E_{1}^{n+1}$ and $M_{1}^{n}$, respectively. Then, the Gauss and Weingarten formulae are given by
\begin{equation}\label{2.1}
\overline\nabla_{X}Y = \nabla_{X}Y + h(X, Y), \hspace{.3 cm} X, Y
\in\Gamma(TM^{n}_{1}),
\end{equation}
\begin{equation}\label{2.2}
\overline\nabla_{X}\xi = -S_{\xi}X, \hspace{.3 cm}  \xi
\in\Gamma(TM^{n}_{1})^{\bot},
\end{equation}
where $h$ is the second fundamental form. It is well known that the
second fundamental form $h$ and shape operator $S$ are related by

\begin{equation}\label{2.3}
\overline{g}(h(X,Y), \xi) = g(S_{\xi}X,Y).
\end{equation}

The mean curvature vector is given by
\begin{equation}\label{2.4}
\vec{H} = \frac{1}{n} \rm trace \emph{h}.
\end{equation}

The Gauss and Codazzi equations are given by
\begin{equation}\label{2.5}
R(X, Y)Z = g(SY, Z) SX - g(SX, Z) SY,
\end{equation}

\begin{equation}\label{2.6}
(\nabla_{X}S)Y = (\nabla_{Y}S)X,
\end{equation}
respectively, where $R$ is the curvature tensor, $S=S_{\xi}$ for
some unit normal vector field $\xi$ and
\begin{equation}\label{2.7}
(\nabla_{X}S)Y = \nabla_{X}(SY)- S(\nabla_{X}Y),
\end{equation}
for all $ X, Y, Z \in \Gamma(TM_{1}^{n})$.

\vspace{.2cm}

   A vector $X$ in $E_{s}^{n+1}$ is called
  spacelike, timelike or lightlike according as \hspace{.2cm}$\overline g(X, X)> 0, \hspace{.2cm} \overline g(X, X)<0$ or \hspace{.2cm} $\overline g(X,
  X)=0,\hspace{.1cm} X\neq 0$, respectively. A non degenerate hypersurface $M^{n}_{r}$ of $E_{s}^{n+1}$ is
  called Riemannian or pseudo-Riemannian according as the induced metric on
  $M_{r}^{n+1}$ from the indefinite metric on $E_{s}^{n+1}$ is definite
  or indefinite. The shape operator of
  pseudo-Riemannian hypersurfaces is not always diagonalizable in contrast to the
  Riemannian hypersurfaces.\vspace{.2cm}

It is known [1, 7] that the matrix representation of the shape
operator of $M^{n}_{1}$ in $E_{1}^{n+1}$ having complex eigenvalues
with respect to a suitable orthonormal base field of the tangent
bundle takes the form

\begin{equation}\label{2.8}
 S = \left(
                            \begin{array}{ccccc}
                               \lambda & -\mu & & \\
                              \mu & \lambda & \\
                              & &     D_{n-2} \\
                             \end{array}
                          \right), \hspace{.5 cm}
 \end{equation}
where $\mu \neq 0$ and
$D_{n-2}=$ diag$\{\lambda_{3},\dots,\lambda_{n}\}$.\\

 The following algebraic lemma will be useful in our study:

\begin{lem}\label{Lemma2.1}\cite[Theorem 4.4, pp. 58--59]{bic6} Let D be a unique factorization domain, and let $f(X) =
a_{0}X^{m} +a_{1}X^{m-1} + � � � + a_{m}, g(X) = b_{0}X^{n} +
b_{1}X^{n-1} + � � � + b_{n}$ be two polynomials in $D[X]$.
Assume that the leading coefficients $a_{0}$ and $b_{0}$ of $f(X)$
and $g(X)$ are not both zero. Then $f(X)$ and $g(X)$ have a non
constant common factor if and only if the resultant $\Re(f, g)$ of
$f$ and $g$ is zero, where
\begin{center}
$\Re(f,g)=
\begin{vmatrix}
  a_{0} & a_{1} & a_{2} & \cdots & a_{m} &   &   &   \\
    & a_{0} & a_{1} & \cdots & \cdots & a_{m} &   &   \\
    &   & \ddots & \ddots & \ddots & \ddots & \ddots &   \\
   &   &   & a_{0} & a_{1} & a_{2} & \cdots & a_{m} \\
  b_{0} & b_{1} & b_{2} & \cdots & b_{n} &   &   &   \\
    & b_{0} & b_{1} & \cdots & \cdots & b_{n} &  &   \\
    &   & \ddots & \ddots & \ddots & \ddots & \ddots &   \\
    &   &   & b_{0} & b_{1} & b_{2} & \cdots & b_{n} \\
\end{vmatrix}
$
\end{center}
where there are n rows of "a" entries and m rows of "b" entries.

\end{lem}

 \section{\textbf{Biconservative Lorentz Hypersurfaces in $E_{1}^{n+1}$ }}

 In this section, we obtain some properties of
biconservative Lorentz hypersurfaces $M_{1}^{n}$ in $E_{1}^{n+1}$
whose shape operator has the form (\ref{2.8}). We assume that $H$ is
not constant which implies that grad$H\neq0$. Hence, there exits an
open connected subset $U$ of $M^{n}_{1}$, with grad$_{p}H \neq 0$
for all $p\in U$. From (\ref{1.2}), it is easy to see that grad$H$
is an eigenvector of the shape operator $S$ with the corresponding
principal curvature $-\frac{n}{2}H$. In view of (\ref{2.8}), the
shape operator $S$ of $M$ satisfies
\begin{multline}\label{3.1}
S(e_{1})=\lambda e_{1}+\mu e_{2},\hspace{.2 cm} S(e_{2})=-\mu
e_{1}+\lambda e_{2},\hspace{.2 cm} S(e_{3})=\lambda_{3} e_{3},\dots,
S(e_{n})=\lambda_{n} e_{n},
 \end{multline}
 with respect to an orthonormal basis $\{e_{1}, e_{2},..., e_{n}\}$ of
$T_{p} M_{1}^{n}$ which satisfies
\begin{equation}\label{3.2}
 g(e_{1},e_{1})=-1, \quad g(e_{i},e_{i})=1, \hspace{.2 cm}
 i=2,3,...,n,
 \end{equation}
and
 \begin{equation}\label{3.3}
  g(e_{i},e_{j})=0, \hspace{.2 cm} for \hspace{.2 cm}  i \neq j.
\end{equation}

We write
\begin{equation}\label{3.4}
    \nabla_{e_{i}}e_{j}=\sum_{k=1}^{n}\omega_{ij}^{k}e_{k},\hspace{2
cm} i, j = 1, 2, ... , n.
\end{equation}

By using (\ref{3.4}) and taking covariant derivatives of (\ref{3.2})
and (\ref{3.3}) with respect to $e_{k}$, we find
     \begin{equation}\label{3.5}
   \begin{split}
    \omega_{ki}^{i}=0,
    \hspace{1 cm}
    \omega_{kj}^{i}=-\omega_{ki}^{j},
    \end{split}
    \end{equation}
for $i \neq j$ and  $i, j, k =1, 2, ... , n.$\\

 In view of (\ref{3.1}), grad$H$ can be chosen in one of
the directions $e_{3},\dots, e_{n}$ and in each direction grad$H$ is
spacelike. Without loss of generality, we may assume $e_{n}$ in the
direction of grad$H$, so $\lambda_{n}= -\frac{nH}{2}.$  We express
grad$H$ as grad$H$ = $-e_{1}(H)e_{1}+e_{2}(H)e_{2}+ \dots
+e_{n}(H)e_{n}$, which gives

\begin{equation}\label{3.6}
e_{n}(H)\neq 0,\quad e_{1}(H)=e_{2}(H)= \dots = e_{n-1}(H)=0.
\end{equation}

Using (\ref{3.4}), (\ref{3.6}) and the fact that $[e_{i} \hspace{.1
cm} e_{j}](H)=0=\nabla_{e_{i}}e_{j}(H)-\nabla_{e_{j}}e_{i}(H),$ for
$i\neq j$ and $i, j \neq n$, we find
\begin{equation}\label{3.7}
\omega_{ij}^{n}=\omega_{ji}^{n}.
\end{equation}

From (\ref{2.7}), (\ref{3.1}), (\ref{3.4}) and (\ref{3.6}), the
Codazzi equation $g((\nabla_{e_{n}}S)e_{a},e_{a}) =
g((\nabla_{e_{a}}S)e_{n},e_{a})$ leads to

\begin{equation}\label{3.8}
e_{n}(\lambda_{a})= (\lambda_{n}-\lambda_{a})\omega_{an}^{a}, \quad
3\leq a\leq n-1.
\end{equation}
Therefore, $\lambda_{n}\neq \lambda_{a}$, otherwise from (\ref{3.8})
we would have $e_{n}(H)=0$ which
contradicts (\ref{3.6}).\\

From $g((\nabla_{X}S)Y,Z) = g((\nabla_{Y}S)X,Z)$, using (\ref{2.7}),
(\ref{3.1}), (\ref{3.4}) and (\ref{3.6}) and the value of
$\lambda_{n}=-\frac{nH}{2}$, we obtain the following equations
showing relations among the
eigenvalues, connection forms and orthonormal frame. These are listed in the last column of Table 1 for $3\leq a,b,c\leq n-1$, $a\neq b\neq c$.\\

\smallskip
 \begin{center}
  Table 1. \ {Evaluation of $g((\nabla _XS)Y, Z)=g((\nabla _YS)X, Z)$ for various values of $e_i$.}
\end{center}

\begin{longtable}{|c | c |c |c |c |}
\hline 
$i$ & X & Y  & Z & Codazzi Equations(T$i$)   \\ [0.5ex] 
\hline\hline 
1 &$e_{1}$ & $e_{2}$ & $e_{1}$ & $e_{2}(\lambda) + e_{1}(\mu)=0 $\\

2 &$e_{1}$ & $e_{2}$ & $e_{2}$ & $e_{1}(\lambda) - e_{2}(\mu)=0 $\\

 3 & $e_{1}$ & $e_{2}$ & $e_{a}$ &
$[\lambda-\lambda_{a}](\omega_{12}^{a}-\omega_{21}^{a})=
\mu(\omega_{22}^{a}+\omega_{11}^{a})$ \\ 

 4 & $e_{1}$ & $e_{2}$ & $e_{n}$ & $[\lambda +
\frac{nH}{2}](\omega_{12}^{n}-\omega_{21}^{n})=
\mu(\omega_{22}^{n}+\omega_{11}^{n})$\\

\hline
5 & $e_{1}$ & $e_{a}$ & $e_{1}$  & $e_{a}(\lambda)= [\lambda_{a}- \lambda]\omega_{1a}^{1} + \mu \omega_{1a}^{2}$\\

 6 & $e_{1}$ & $e_{a}$ & $e_{2}$ & $e_{a}(\mu)= [\lambda_{a}-
\lambda]\omega_{1a}^{2} - \mu \omega_{1a}^{1}.
$ \\

 7 & $e_{1}$ & $e_{a}$ & $e_{a}$ & $ e_{1}(\lambda_{a})= [\lambda-
\lambda_{a}]\omega_{a1}^{a} + \mu \omega_{a2}^{a}
$ \\

 8 & $e_{1}$ & $e_{a}$ & $e_{n}$ & $[\lambda_{a} +
\frac{nH}{2}]\omega_{1a}^{n} = [\lambda +
\frac{nH}{2}]\omega_{a1}^{n}
+ \mu \omega_{a2}^{n}$\\

 9 & $e_{1}$ & $e_{a}$ & $e_{b}$ & $[\lambda_{a}-
\lambda_{b}]\omega_{1a}^{b} = [\lambda - \lambda_{b}]\omega_{a1}^{b}
+ \mu \omega_{a2}^{b}$ \\

\hline
10 & $e_{1}$ & $e_{n}$ & $e_{1}$ & $ -(\lambda + \frac{nH}{2})\omega_{1n}^{1}+ \mu \omega_{1n}^{2}=e_{n}(\lambda)$\\

 11 & $e_{1}$ & $e_{n}$ & $e_{2}$ & $-(\lambda +
\frac{nH}{2})\omega_{1n}^{2}- \mu
\omega_{1n}^{1}=e_{n}(\mu)$\\

12 & $e_{1}$ & $e_{n}$ & $e_{n}$ & $ (\lambda + \frac{nH}{2})\omega_{n1}^{n}+ \mu \omega_{n2}^{n}=0$\\

\hline
13 & $e_{2}$ & $e_{a}$ & $e_{1}$ & $ -e_{a}(\mu)= [\lambda_{a}-\lambda]\omega_{2a}^{1}+ \mu \omega^{2}_{2a}$\\

14 & $e_{2}$ & $e_{a}$ & $e_{2}$ & $ e_{a}(\lambda)=[\lambda_{a}-\lambda]\omega_{2a}^{2}- \mu \omega^{1}_{2a}$\\

15 & $e_{2}$ & $e_{a}$ & $e_{a}$ & $ e_{2}(\lambda_{a})= [\lambda-\lambda_{a}]\omega_{a2}^{a}- \mu \omega^{a}_{a1}$\\

 16 & $e_{2}$ & $e_{a}$ & $e_{n}$ &
$[\lambda_{a}+\frac{nH}{2}]\omega_{2a}^{n} = [\lambda +
\frac{nH}{2}]\omega_{a2}^{n}
- \mu \omega_{a1}^{n}$\\

17 & $e_{2}$ & $e_{a}$ & $e_{b}$ & $(\lambda_{a}-\lambda_{b})\omega_{2a}^{b}= [\lambda-\lambda_{b}]\omega_{a2}^{b}- \mu \omega^{b}_{a1}$\\

\hline
18 & $e_{2}$ & $e_{n}$ & $e_{1}$ & $ -(\lambda + \frac{nH}{2})\omega_{2n}^{1}+ \mu \omega_{2n}^{2}= -e_{n}(\mu)$\\

 19 & $e_{2}$ & $e_{n}$ & $e_{2}$ & $-(\lambda +
\frac{nH}{2})\omega_{2n}^{2}- \mu
\omega_{2n}^{1}=e_{n}(\lambda)$\\

20 & $e_{2}$ & $e_{n}$ & $e_{n}$ & $ (\lambda + \frac{nH}{2})\omega_{n2}^{n}- \mu \omega_{n1}^{n}=0$\\

\hline 21 & $e_{a}$ & $e_{b}$ & $e_{1}$ & $(\lambda_{b} -
\lambda)\omega_{ab}^{1}+ \mu \omega_{ab}^{2}=
[\lambda_{a}-\lambda]\omega_{ba}^{1}+ \mu \omega_{ba}^{2}$\\

 22 & $e_{a}$ & $e_{b}$ & $e_{2}$ & $(\lambda_{b} -
\lambda)\omega_{ab}^{2}- \mu \omega_{ab}^{1}=
[\lambda_{a}-\lambda]\omega_{ba}^{2}- \mu \omega_{ba}^{1}$\\

 23 & $e_{a}$ & $e_{b}$ & $e_{n}$ &
$[\lambda_{b}+\frac{nH}{2}]\omega_{ab}^{n}=
[\lambda_{a}+\frac{nH}{2}]\omega_{ba}^{n}$\\

 24 & $e_{a}$ & $e_{b}$ & $e_{c}$ &
$[\lambda_{b}-\lambda_{c}]\omega_{ab}^{c}=
[\lambda_{a}-\lambda_{c}]\omega_{ba}^{c}$\\

\hline 25 & $e_{a}$ & $e_{n}$ & $e_{1}$ & $-(\lambda +
\frac{nH}{2})\omega_{an}^{1}+ \mu \omega_{an}^{2}=
[\lambda_{a}-\lambda]\omega_{na}^{1}+ \mu \omega_{na}^{2}$\\

 26 & $e_{a}$ & $e_{n}$ & $e_{2}$ & $-(\lambda +
\frac{nH}{2})\omega_{an}^{2}- \mu \omega_{an}^{1}=
[\lambda_{a}-\lambda]\omega_{na}^{2}- \mu \omega_{na}^{1}$\\

27 & $e_{a}$ & $e_{n}$ & $e_{a}$ & $e_{n}(\lambda_{a})= -[\frac{nH}{2}+ \lambda_{a}]\omega_{an}^{a}$\\

28 & $e_{a}$ & $e_{n}$ & $e_{n}$ & $ \omega_{na}^{n}=0$\\

 29 & $e_{a}$ & $e_{n}$ & $e_{b}$ &
$-[\lambda_{b}+\frac{nH}{2}]\omega_{an}^{b}=
[\lambda_{a}-\lambda_{b}]\omega_{na}^{b}$\\[1ex] 
\hline 
\end{longtable}

\medskip

 Using T23, T29, (\ref{3.7}) and (\ref{3.5}), we have
\begin{equation}\label{3.9}
         \omega_{ab}^{n}=\omega_{ba}^{n}=\omega_{an}^{b}=\omega_{nb}^{a}=0
         \quad \mbox{for} \quad \lambda_{a}\neq \lambda_{b}.
\end{equation}

Note that all the connection coefficients vanish in (\ref{3.9}) for
$\lambda_{a}\neq \lambda_{b}$, will also vanish for $\lambda_{a}=
\lambda_{b}$ except $\omega_{nb}^{a}$. Now, equating T10, T19 and
T11, T18 and using (\ref{3.5}), we find
  \begin{equation}\label{3.10}
         \omega_{22}^{n}=\omega_{11}^{n}, \quad \omega_{12}^{n}= -\omega_{21}^{n},
\end{equation}
which by use of T4, (\ref{3.7}) and (\ref{3.5}) gives

 \begin{equation}\label{3.11}
         \omega_{22}^{n}=\omega_{11}^{n}= \omega_{12}^{n}= \omega_{21}^{n}=\omega_{2n}^{1}=
         \omega_{1n}^{2}=\omega_{2n}^{2}=
         \omega_{1n}^{1}=0.
\end{equation}

Similarly, using T5, T14, T6, T13 and (\ref{3.5}), we find
\begin{equation}\label{3.12}
         \omega_{22}^{a}=\omega_{11}^{a}, \quad \omega_{12}^{a}= -\omega_{21}^{a}.
\end{equation}

 Using T12, T20, T28 and (\ref{3.5}), we get
\begin{equation}\label{3.13}
         \omega_{n1}^{n}=\omega_{n2}^{n}= \omega_{nn}^{1}=
         \omega_{nn}^{2}= \omega_{nn}^{a}=0.
\end{equation}

Solving T8, T16 by using (\ref{3.7}) and (\ref{3.5}), we obtain
\begin{equation}\label{3.14}
         \omega_{1a}^{n}=\omega_{2a}^{n}=\omega_{a1}^{n}=
         \omega_{a2}^{n}= \omega_{1n}^{a}=\omega_{2n}^{a}=\omega_{an}^{1}=
         \omega_{an}^{2}=0.\\
\end{equation}

Now, solving T25, T26 by using (\ref{3.14}) and (\ref{3.5}), we
obtain
\begin{equation}\label{3.15}
         \omega_{na}^{1}=\omega_{na}^{2}=\omega_{n1}^{a}=
         \omega_{n2}^{a}=0.
\end{equation}

Replacing $e_{a}$ and $e_{b}$ in T9 and T17 and using T9, T17 and
(\ref{3.5}), we obtain for $\lambda_{a}\neq \lambda_{b}$
\begin{equation}\label{3.16}
[\lambda - \lambda_{a}]\omega_{ba}^{1} + \mu \omega_{ba}^{2} =
[\lambda - \lambda_{b}]\omega_{ab}^{1} + \mu \omega_{ab}^{2},
\end{equation}

\begin{equation}\label{3.17}
[\lambda-\lambda_{a}]\omega_{ba}^{2}- \mu \omega^{1}_{ba}=
[\lambda-\lambda_{b}]\omega_{ab}^{2}- \mu \omega^{1}_{ab}
\end{equation}

Using (\ref{3.16}) and T21, we get

\begin{equation}\label{3.18}
         (\lambda_{a}-\lambda)\omega_{ba}^{1}=(\lambda_{b}-\lambda)\omega_{ab}^{1}, \quad \omega_{ab}^{2}=
         \omega_{ba}^{2}.
\end{equation}

Similarly, Using (\ref{3.17}) and T22, we get

\begin{equation}\label{3.19}
         (\lambda_{a}-\lambda)\omega_{ba}^{2}=(\lambda_{b}-\lambda)\omega_{ab}^{2}, \quad \omega_{ab}^{1}=
         \omega_{ba}^{1}.
\end{equation}

Combining (\ref{3.18}) and (\ref{3.19}), T9, T17 and using
(\ref{3.5}), we obtain

\begin{equation}\label{3.20}
         \omega_{ab}^{1}= \omega_{ab}^{2}=
         \omega_{a1}^{b}= \omega_{a2}^{b}=
         \omega_{1a}^{b}=\omega_{2a}^{b}=0 \quad \mbox{for} \quad \lambda_{a}\neq \lambda_{b}.
\end{equation}

Also, it can be easily seen that all the connection coefficients
vanish in (\ref{3.20}) for $\lambda_{a}\neq \lambda_{b}$, will also
vanish for $\lambda_{a}= \lambda_{b}$ except $\omega_{1a}^{b}$ and
$\omega_{2a}^{b}$.

\vspace{.5 cm}

By using the above relations, we obtain the following:
\begin{lem}\label{Lemma3.1}
 Let $M^{n}_{1}$ be a biconservative hypersurface in $E^{n+1}_{1}$, whose shape
 operator has the form (\ref{2.8}) with respect to suitable orthonormal basis  $\{e_{1},
e_{2}, ..., e_{n}\}$. If grad$H$ in the direction of $e_{n }$, then

\begin{center}
      $$\nabla_{e_{1}}e_{1}=\sum_{p\neq 1,n}\omega_{11}^{p}e_{p},\hspace{.2 cm}
      \nabla_{e_{1}}e_{2}=\sum_{p\neq 2,n}\omega_{12}^{p}e_{p},\hspace{.2 cm}
      \nabla_{e_{1}}e_{n}=0,\hspace{.2 cm}
      \nabla_{e_{2}}e_{1}= \sum_{p\neq 1,n}\omega_{21}^{p}e_{p},$$\hspace{.2 cm}
      $$\nabla_{e_{2}}e_{2}=\sum_{p\neq 2,n}\omega_{22}^{p}e_{p},\hspace{.2 cm}
      \nabla_{e_{2}}e_{n}=0,\hspace{.2 cm}
      \nabla_{e_{a}}e_{1}=\sum_{p\neq 1,b,n}\omega_{a1}^{p}e_{p},\hspace{.2 cm}
      \nabla_{e_{a}}e_{2}=\sum_{p\neq 2,b,n}\omega_{a2}^{p}e_{p},$$\hspace{.2 cm}
      $$\nabla_{e_{a}}e_{a}=\sum_{p\neq a}\omega_{aa}^{p}e_{p},\hspace{.2 cm}
      \nabla_{e_{a}}e_{b}=\sum_{p\neq 1,2,b,n}\omega_{ab}^{p}e_{p},\hspace{.2 cm}
      \nabla_{e_{a}}e_{n}= \sum_{p\neq 1,2,b,n}\omega_{an}^{p}e_{p},$$\hspace{.2 cm}
      $$\nabla_{e_{n}}e_{1}=\omega_{n1}^{2}e_{2},\hspace{.2 cm}
      \nabla_{e_{n}}e_{2}=\omega_{n2}^{1}e_{1},\hspace{.2 cm}
      \nabla_{e_{n}}e_{n}=0,$$\hspace{.1 cm}
      $$
  \nabla_{e_{1}}e_{a}=\begin{cases}
               \omega_{1a}^{1}e_{1}+\omega_{1a}^{2}e_{2}, \hspace{2.35 cm} \lambda_{a}\neq \lambda_{b}, \\
               \omega_{1a}^{1}e_{1}+\omega_{1a}^{2}e_{2}+\omega_{1a}^{b}e_{b}, \hspace{1.0 cm} \lambda_{a}= \lambda_{b},
            \end{cases}
$$
 $$
  \nabla_{e_{2}}e_{a}=\begin{cases}
               \omega_{2a}^{1}e_{1}+\omega_{2a}^{2}e_{2}, \hspace{2.35 cm} \lambda_{a}\neq \lambda_{b}, \\
               \omega_{2a}^{1}e_{1}+\omega_{2a}^{2}e_{2}+\omega_{2a}^{b}e_{b}, \hspace{1.0 cm} \lambda_{a}= \lambda_{b},
            \end{cases}
$$
$$
  \nabla_{e_{n}}e_{a}=\begin{cases}
               0, \hspace{1.7 cm} \lambda_{a}\neq \lambda_{b}, \\
               \omega_{na}^{b}e_{b}, \hspace{1.0 cm} \lambda_{a}=
               \lambda_{b}.
            \end{cases}
$$
\end{center}
\end{lem}

\vspace{.2 cm}

 Using Lemma \ref{Lemma3.1}, (\ref{2.5}), (\ref{3.1}), we evaluate $g(R(e_{2},e_{n})e_{2},e_{n})$,
 $g(R(e_{a},e_{n})e_{a},e_{n})$,\\ $g(R(e_{a},e_{l})e_{a},e_{n})$ and we obtain the following:
\begin{lem}\label{Lemma3.2}
Let $M^{n}_{1}$ be a biconservative hypersurface in $E^{n+1}_{1}$,
whose shape operator has the form (\ref{2.8}) with respect to
suitable orthonormal basis  $\{e_{1}, e_{2}, ..., e_{n}\}$. If
grad$H$ in the direction of $e_{n }$, then we have
\begin{equation}\label{3.21}
 \lambda=0,
\end{equation}
\begin{equation}\label{3.22}
 e_{n}(\omega_{aa}^{n})- (\omega_{aa}^{n})^{2} = -\frac{nH}{2}\lambda_{a}, \hspace{.2 cm}a =
 3, 4, \dots n-1,
\end{equation}
 \begin{equation}\label{3.23}
 e_{l}(\omega_{aa}^{n})= \omega_{aa}^{l} \omega_{aa}^{n}, \hspace{.2 cm}l = 1,
 2, \quad a=3,4, \dots n-1.
\end{equation}
\end{lem}
\emph{\textbf{Proof:}} To prove equations
(\ref{3.21})$\sim$(\ref{3.23}), we consider the following cases:\\

\textbf{Case A:} If shape operator (\ref{2.8}) has all distinct
principal curvatures \emph{i.e} $\lambda_{a}\neq \lambda_{b}$ for
all $a,b=3,4, \dots, n-1$.

 From (\ref{2.5}), (\ref{3.1}) and using Lemma
\ref{Lemma3.1}, we evaluate
\begin{center} $\bullet$ $g(R(e_{2},e_{n})e_{2},e_{n})= g(Ae_{n},
e_{2}) g(Ae_{2}, e_{n}) - g(Ae_{2}, e_{2})
g(Ae_{n},e_{n})=\frac{nH}{2}\lambda,$
 \end{center}
 or,
\begin{align*}
 \frac{nH}{2}\lambda &=
g\left(\nabla_{e_{2}}\nabla_{e_{n}}e_{2}-\nabla_{e_{n}}\nabla_{e_{2}}e_{2}-\nabla_{[e_{2}
e_{n}]}e_{2},
e_{n}\right)\\
&=g\left(\nabla_{e_{2}}(\omega_{n2}^{1}e_{1})-\nabla_{e_{n}}\bigg(\sum_{p\neq
2,n}\omega_{22}^{p}e_{p}\bigg)-\nabla_{(-\omega_{n2}^{1}e_{1})}e_{2}, e_{n}\right)\\
&=
g\Bigg(e_{2}\left(\omega_{n2}^{1}\right)e_{1}+\omega_{n2}^{1}\sum_{p\neq
1,n}\omega_{21}^{p}e_{p}-\sum_{p\neq
2,n}\left(e_{n}\left(\omega_{22}^{p}\right)e_{p}+\omega_{22}^{p}\nabla_{e_{n}}e_{p}\right)\\
&\qquad+\omega_{n2}^{1}\sum_{p\neq 2,n}\omega_{12}^{p}e_{p},
e_{n}\Bigg)=0.\end{align*}
\begin{center}
$ \bullet \ g(R(e_{a},e_{n})e_{a},e_{n})= g(Ae_{n}, e_{a}) g(Ae_{a},
e_{n})- g(Ae_{a}, e_{a}) g(Ae_{n}, e_{n})=\frac{nH}{2}\lambda_{a},$
\end{center}
  hence,\begin{align*}\frac{nH}{2}\lambda_{a}&=
g\left(\nabla_{e_{a}}\nabla_{e_{n}}e_{a}-\nabla_{e_{n}}\nabla_{e_{a}}e_{a}-\nabla_{[e_{a}
e_{n}]}e_{a},
e_{n}\right)\\&=g\left(-\nabla_{e_{n}}\bigg(\sum_{p\neq a,
p=1}^{n}\omega_{aa}^{p}e_{p}\bigg)-\nabla_{(\omega_{an}^{a}e_{a})}e_{a},
e_{n}\right)\\&=g\left(-\sum_{p\neq
a,p=1}^{n}\left(e_{n}(\omega_{aa}^{p})e_{p}+\omega_{aa}^{p}\nabla_{e_{n}}e_{p}\right)-\omega_{an}^{a}\sum_{p\neq
a,p=1}^{n}\omega_{aa}^{p}e_{p}, e_{n}\right)\\&=
-e_{n}(\omega_{aa}^{n})-\omega_{an}^{a}\omega_{aa}^{n}.\end{align*}
\begin{center}
$\bullet\ g(R(e_{a},e_{1})e_{a},e_{n})= g(Ae_{1}, e_{a}) g(Ae_{a},
e_{n}) - g(Ae_{a}, e_{a}) g(Ae_{1}, e_{n})=0,$
\end{center}
 hence,
 \begin{align*}0&=
g\left(\nabla_{e_{a}}\nabla_{e_{1}}e_{a}-\nabla_{e_{1}}\nabla_{e_{a}}e_{a}-\nabla_{[e_{a}
e_{1}]}e_{a}, e_{n}\right)\\&=g\left(\nabla_{e_{a}}\bigg(\sum_{p\neq
a,b,n}\omega_{1a}^{p}e_{p}\bigg)-\nabla_{e_{1}}\bigg(\sum_{p\neq
a}\omega_{aa}^{p}e_{p}\bigg)-\nabla_{\nabla_{e_{a}}
e_{1}-\nabla_{e_{1}}
e_{a}}e_{a}, e_{n}\right)\\
&=g\left(\sum_{p\neq a,
b,n}(e_{a}(\omega_{1a}^{p})e_{p}+\omega_{1a}^{p}\nabla_{e_{a}}e_{p}),
e_{n}\right)-g\left(\sum_{p\neq
a}(e_{1}(\omega_{aa}^{p})e_{p}+\omega_{aa}^{p}\nabla_{e_{1}}e_{p}),
e_{n}\right)\\
&\qquad-\displaystyle{g\bigg(\nabla_{\nabla_{e_{a}}
e_{1}-\nabla_{e_{1}} e_{a}}e_{a}, e_{n}\bigg)}\end{align*}

Now, it is $$ g\bigg(\sum_{p\neq a,
b,n}(e_{a}(\omega_{1a}^{p})e_{p}+\omega_{1a}^{p}\nabla_{e_{a}}e_{p}),
e_{n}\bigg)=0,$$ as $\nabla_{e_{a}}e_{p}$ does not have a component
along $e_{n}$ for $p \neq a,b,n$. Similarly,
$$g\bigg(\sum_{p\neq
a}(e_{1}(\omega_{aa}^{p})e_{p}+\omega_{aa}^{p}\nabla_{e_{1}}e_{p}),
e_{n}\bigg)= e_{1}(\omega_{aa}^{n})$$ and
$$g\bigg(\nabla_{\nabla_{e_{a}} e_{1}}e_{a},
e_{n}\bigg)=g\bigg(\sum_{p\neq 1,b,n}\omega_{a1}^{p}\nabla_{e_{p}}
e_{a}, e_{n}\bigg)=-\omega_{a1}^{a}\omega_{aa}^{n},$$
\\$$ g\bigg(\nabla_{\nabla_{e_{1}} e_{a}}e_{a}, e_{n}\bigg)=g\bigg(\sum_{ p\neq
a,b,n}\omega_{1a}^{p}\nabla_{e_{p}} e_{a}, e_{n}\bigg)=0.$$

Therefore,
$$e_{1}(\omega_{aa}^{n})-\omega_{aa}^{1}\omega_{aa}^{n}=0.$$

In the same way, by evaluating $g(R(e_{a},e_{2})e_{a},e_{n})$, we
obtain
$$e_{2}(\omega_{aa}^{n})-\omega_{aa}^{2}\omega_{aa}^{n}=0.$$
\\
\textbf{Case B:} If we have $\lambda_{a}=\lambda_{b}$ for any
$a,b\in \{3,4,\dots,n-1\}$.

Also, in this case we can prove easily
(\ref{3.21})$\sim$(\ref{3.23}) in the same way as we have proved in
Case A.

Hence, the proof of Lemma is completed.\\
\\

Now, using Lemma \ref{Lemma3.1}, Lemma \ref{Lemma3.2}, Table 1 and
(\ref{2.5}), we can obtain the following:

\begin{lem}\label{Lemma3.3}
 Let $M^{n}_{1}$ be a biconservative hypersurface in $E^{n+1}_{1}$, whose shape operator has the form (\ref{2.8})
  with respect to a suitable orthonormal basis $\{e_{1},
e_{2}, ..., e_{n}\}$. If grad$H$ in the direction of $e_{n}$, then
$\mu$ is constant.
\end{lem} \vspace{.2 cm}

\textbf{Proof:} We use T3, T5, (\ref{3.12}), (\ref{3.21}) and
(\ref{3.5}) , we obtain that

\begin{equation}\label{3.24}
         \omega_{12}^{a}(\lambda_{a}^{2}-
         \mu^{2})=0.
\end{equation}

If $\lambda_{a}^{2}= \mu^{2}$ for any $a=3,4, \dots n-1$, then
$e_{n}(\lambda_{a})=0$ as $e_{n}(\mu)=0$ from T11 and (\ref{3.11}).
Therefore from T27, we have $\omega^{n}_{aa}=0$ which by using
(\ref{3.22}) implies that $\lambda_{a}=0$ as $H\neq 0.$ This
contradicts to the fact that $\mu$ is non zero. So, from
(\ref{3.24}), T3, (\ref{3.12}) and (\ref{3.5}) we get
\begin{equation}\label{3.25}
\omega_{21}^{a}=\omega_{22}^{a}=\omega_{11}^{a}=\omega_{2a}^{1}=\omega_{2a}^{2}=\omega_{1a}^{1}=\omega^{a}_{12}=\omega^{2}_{1a}=0,\quad
a=3,4,\dots n-1.
\end{equation}

Now, from T1, T2, T6, T18, (\ref{3.11}), (\ref{3.21}) and
(\ref{3.25}), it follows that $\mu$ is constant in every direction
which completes the
proof.\\

In a similar way we use Lemma \ref{Lemma3.1}, (\ref{2.5}),
(\ref{3.1}), (\ref{3.25}), to evaluate
$g(R(e_{a},e_{2})e_{a},e_{1})$, $g(R(e_{a},e_{1})e_{a},e_{2})$ and
we obtain the following:
\begin{lem}\label{Lemma3.4}
Let $M^{n}_{1}$ be biconservative hypersurface in $E^{n+1}_{1}$,
whose shape operator has the form (\ref{2.8})
 with respect to a suitable orthonormal basis. If grad$H$ in the direction
of $e_{n }$, then the following relations are valid:
\begin{equation}\label{3.26}
 e_{2}(\omega_{aa}^{1})+ \omega_{aa}^{2}(\omega_{22}^{1}-\omega_{aa}^{1}) = -\mu \lambda_{a},
\end{equation}
 \begin{equation}\label{3.27}
 e_{1}(\omega_{aa}^{2})+ \omega_{aa}^{1}(\omega_{11}^{2}-\omega_{aa}^{2}) = -\mu \lambda_{a},
\end{equation}
for $a=3,4, \dots n-1.$
\end{lem}
\emph{\textbf{Proof:}} This is similar to the proof of Lemma
\ref{Lemma3.2}.

\section{\textbf{Biconservative Lorentz hypersurfaces with at most five distinct principal curvature}}

Now, we are in the position to prove our main Theorem
\ref{Theorem1.1}. Therefore, in this section, we study
biconservative Lorentz hypersurface in $E_{1}^{n+1}$ having shape
operator (\ref{2.8}) with at
most five distinct principal curvatures.\\

$$\textbf{Proof of Theorem 1.1}$$

\emph{Case-(i): Five distinct principal curvatures}

 We use (\ref{3.21}) and see easily that the eigenvalues of
 shape operator (\ref{2.8}) are $ \pm \sqrt{-1} \mu,$ $
 \lambda_{3}, \dots ,\lambda_{n}$. So, under the assumption that shape operator
 (\ref{2.8}) has five distinct eigenvalues for $n\geq 5$, we may assume that $\lambda_{3}=\lambda_{4}= \ldots =\lambda_{r}$ and $\lambda_{r+1}=\lambda_{r+2}= \ldots
 =\lambda_{n-1}$.

 So, using (\ref{2.4}), we have trace$S=nH$ and using (\ref{3.21}) and the value
of $\lambda_{n}=-\frac{nH}{2}$, we get

 \begin{equation}\label{4.1}
 (r-2)\lambda_{3}+(n-r-1)\lambda_{n-1}= \frac{3nH}{2}.
 \end{equation}

 Expressions (\ref{3.1}) reduce to
\begin{equation}\label{4.2}
\begin{array}{rcl}
S(e_{1})=\mu e_{2},\hspace{.2 cm} S(e_{2})=-\mu e_{1},\hspace{.2 cm}
S(e_{A})=\lambda_{3} e_{A},\\\hspace{.2 cm} S(e_{B})=\lambda_{n-1}
e_{B}, \hspace{.2 cm}, S(e_{n})=\lambda_{n} e_{n},
\end{array}
 \end{equation}
where $A=3,4, \dots ,r$, $B= r+1, r+2, \dots ,n-1$.

Differentiating (\ref{4.1}) along $e_{n}$ and using (\ref{4.1}), T27
for $a=A,B$, we get
\begin{equation}\label{4.3}
3ne_{n}(H)=[n(n-r+2)H-2(r-2)\lambda_{3}]\omega^{n}_{BB}
+(r-2)(2\lambda_{3}+nH)\omega^{n}_{AA}.
\end{equation}

Using Lemma \ref{Lemma3.1}, (\ref{3.6}) and the fact that
$[e_{i}\hspace{.1 cm}
e_{n}](H)=0=\nabla_{e_{i}}e_{n}(H)-\nabla_{e_{n}}e_{i}(H)$ for $i
=1, 2, \ldots, n-1$, we obtain
      \begin{equation}\label{4.4}
      e_{i}e_{n}(H)= 0.
      \end{equation}

  Differentiating (\ref{4.1}) along $e_{1}$ and using T7 for $a=A, B$  and (\ref{4.1}), we obtain

\begin{equation}\label{4.5}
         2(r-2)[\lambda_{3}\omega_{AA}^{1}-\mu\omega_{AA}^{2}]+[(3nH-2(r-2)\lambda_{3})\omega_{BB}^{1}-2\mu(n-r-1)\omega_{BB}^{2}]=0.
       \end{equation}

      Similarly, differentiating (\ref{4.1}) along $e_{2}$ and using T15 for $a=A, B$  and (\ref{4.1}), we obtain that

\begin{equation}\label{4.6}
         2(r-2)[\lambda_{3}\omega_{AA}^{2}+\mu\omega_{AA}^{1}]+[(3nH-2(r-2)\lambda_{3})\omega_{BB}^{2}+2\mu(n-r-1)\omega_{BB}^{1}]=0.
       \end{equation}

       Multiplying (\ref{4.5}) and (\ref{4.6}) by $\lambda_{3}$ and $\mu$
       respectively and then adding, we get

       \begin{multline}\label{4.7}
         2(r-2)(\lambda_{3}^{2}+\mu^{2})\omega_{AA}^{1}+[\{\lambda_{3}(3nH-2(r-2)\lambda_{3})+2\mu^{2}(n-r-1)\}\omega_{BB}^{1}\\
         +\{\mu(3nH-2(r-2)\lambda_{3})-2\mu \lambda_{3}(n-r-1)\}\omega_{BB}^{2}]=0.
       \end{multline}

        Now, multiplying (\ref{4.5}) and (\ref{4.6}) by $\mu$ and $\lambda_{3}$
       respectively and subtracting, we get

       \begin{multline}\label{4.8}
         2(r-2)(\lambda_{3}^{2}+\mu^{2})\omega_{AA}^{2}+[\{\lambda_{3}(3nH-2(r-2)\lambda_{3})+2\mu^{2}(n-r-1)\}\omega_{BB}^{2}\\
         +\{-\mu(3nH-2(r-2)\lambda_{3})+2\mu \lambda_{3}(n-r-1)\}\omega_{BB}^{1}]=0.
       \end{multline}

       Differentiating (\ref{4.3}) along $e_{1}$ and using (\ref{4.4}), (\ref{4.1}) and (\ref{3.23}) for $a=A, B$ and $l=1$, we obtain

       \begin{multline}\label{4.9}
         2(r-2)[\lambda_{3}\omega_{AA}^{1}-\mu\omega_{AA}^{2}](\omega^{n}_{AA}-\omega^{n}_{BB})+(n(n-r+2)H-2(r-2)\lambda_{3})\omega_{BB}^{n}\omega_{BB}^{1}\\
         +(r-2)(2\lambda_{3}+nH)\omega_{AA}^{n}\omega_{AA}^{1}=0.
       \end{multline}

      Similarly, differentiating (\ref{4.3}) along $e_{2}$ and using (\ref{4.4}), (\ref{4.1}) and (\ref{3.23}) for $a=A, B$ and $l=2$ we obtain

       \begin{multline}\label{4.10}
         2(r-2)[\lambda_{3}\omega_{AA}^{2}+\mu\omega_{AA}^{1}](\omega^{n}_{AA}-\omega^{n}_{BB})+(n(n-r+2)H-2(r-2)\lambda_{3})\omega_{BB}^{n}\omega_{BB}^{2}\\
         +(r-2)(2\lambda_{3}+ nH)\omega_{AA}^{n}\omega_{AA}^{2}=0.
       \end{multline}

       Eliminating $\omega_{AA}^{1}$ and $\omega_{AA}^{2}$ from
       (\ref{4.9}) using (\ref{4.5}) and (\ref{4.7}), we obtain

       \begin{multline}\label{4.11}
        \omega_{AA}^{n}[2\mu(n-r-1)\omega^{2}_{BB}-(3nH-2(r-2)\lambda_{3})\omega_{BB}^{1}-\frac{(2\lambda_{3}+nH)}{2(\lambda_{3}^{2}+\mu^{2})}
        (\{\lambda_{3}(3nH-2(r-2)\lambda_{3})+\\2\mu^{2}(n-r-1)\}\omega_{BB}^{1}
         +\{\mu(3nH-2(r-2)\lambda_{3})-2\mu \lambda_{3}(n-r-1)\}\omega_{BB}^{2})]
         +\omega_{BB}^{n}[-2\mu(n-r-1)\omega^{2}_{BB}\\+(n(n-r+5)H-4(r-2)\lambda_{3})\omega_{BB}^{1}]=0.
       \end{multline}

       Similarly, eliminating $\omega_{AA}^{1}$ and $\omega_{AA}^{2}$ from
       (\ref{4.10}) using (\ref{4.6}) and (\ref{4.8}), we obtain

\begin{multline}\label{4.12}
        \omega_{AA}^{n}[-2\mu(n-r-1)\omega^{1}_{BB}-(3nH-2(r-2)\lambda_{3})\omega_{BB}^{2}-\frac{(2\lambda_{3}+nH)}{2(\lambda_{3}^{2}+\mu^{2})}
        \huge{ ( }\{\lambda_{3}(3nH-2(r-2)\lambda_{3})+\\2\mu^{2}(n-r-1)\}\omega_{BB}^{2}
         +\{-\mu(3nH-2(r-2)\lambda_{3})+2\mu \lambda_{3}(n-r-1)\}\omega_{BB}^{1})]
         +\omega_{BB}^{n}[2\mu(n-r-1)\omega^{1}_{BB}\\+(n(n-r+5)H-4(r-2)\lambda_{3})\omega_{BB}^{2}]=0.
       \end{multline}

       Now, eliminating $\omega^{n}_{AA}$ and $\omega^{n}_{BB}$ from
       (\ref{4.11}) and (\ref{4.12}), we get

       \begin{multline}\label{4.13}
        [(\omega_{BB}^{1})^{2}+(\omega_{BB}^{2})^{2}][2PQ(\lambda_{3}^{2}+\mu^{2})+Q(2\lambda_{3}+nH)(\lambda_{3}P-\mu
        R)-2PR(\lambda_{3}^{2}+\mu^{2})\\
        -P(2\lambda_{3}+nH)(\lambda_{3}R+ \mu P)]=0,
       \end{multline}
       where $P=2\mu(n-r-1)$, $Q=n(n-r+5)H-4(r-2)\lambda_{3}$ and
       $R=3nH-2(r-2)\lambda_{3}.$\\

       We claim $(\omega_{BB}^{1})^{2}+(\omega_{BB}^{2})^{2}\neq
       0$.\\

      Indeed, if $(\omega_{BB}^{1})^{2}+(\omega_{BB}^{2})^{2}=0$, we have,
      $\omega_{BB}^{1}=\omega_{BB}^{2}=0$ as connection coefficients are real numbers. Then, using (\ref{4.7}) and (\ref{4.8}),
      we have $\omega_{AA}^{1}=\omega_{AA}^{2}=0$.

     Therefore, using Lemma \ref{Lemma3.4} for $a=A, B$ , we obtain

      \begin{equation}\label{4.14}
         \lambda_{3}\mu=0, \quad
         \lambda_{n-1}\mu=0,
       \end{equation}
       respectively which implies $\lambda_{3}=\lambda_{n-1}=0$.
       Using T27 for $a=A, B$, we obtain that
       $\omega_{AA}^{n}=\omega_{BB}^{n}=0$. Also, from (\ref{4.3}) we have $e_{n}(H)=0$ which is a contradiction. Hence our claim is proved.\\

Therefore, we have
\begin{multline}\label{4.15}
        f(\lambda_{3},H)\equiv2PQ(\lambda_{3}^{2}+\mu^{2})+Q(2\lambda_{3}+nH)(\lambda_{3}P-\mu
        R)-2PR(\lambda_{3}^{2}+\mu^{2})
        -P(2\lambda_{3}+nH)(\lambda_{3}R+ \mu P)=0.
       \end{multline}

              Differentiating (\ref{4.15}) along $e_{1}$ and $e_{2}$ and using Lemma \ref{Lemma3.3}, we have

       \begin{equation}\label{4.16}
         e_{1}(\lambda_{3})g(\lambda_{3},H)=0,
       \end{equation}
       and
\begin{equation}\label{4.17}
         e_{2}(\lambda_{3})g(\lambda_{3},H)=0,
       \end{equation}
       respectively, where $g(\lambda_{3},H)= 4P\lambda_{3}(Q-R)-4P(\lambda_{3}^{2}+\mu^{2})(r-2)+2(PQ\lambda_{3}-QR\mu-\lambda_{3}PR-P^{2}\mu)
       +(2\lambda_{3}+nH)(PQ-2\lambda_{3}P(r-2)+2(r-2)(2R+Q)\mu-PR)$.

Now, if $g(\lambda_{3},H)\neq 0$, we have $e_{1}(\lambda_{3})=0$ and
$e_{2}(\lambda_{3})=0$ which implies from T7, T15 for $a=A$,
(\ref{4.5}) and (\ref{4.6}) that
$\omega_{BB}^{1}=\omega_{BB}^{2}=\omega_{AA}^{1}=\omega_{AA}^{2}=0$.
As we have already proved from (\ref{4.14}) that it arises to a
contradiction.

Therefore, we have

\begin{equation}\label{4.18}
         g(\lambda_{3},H)= 0,
       \end{equation}
  which is a polynomial equation in $\lambda_{3}$ and $H$.

We rewrite $f(\lambda_{3},H)$, $g(\lambda_{3},H)$ as polynomials
$f_{H}(\lambda_{3}), g_{H}(\lambda_{3})$ of $\lambda_{3}$ with
coefficients in the polynomial ring $R[H]$ over $\mathbb{R}$. Since
$f_{H}(\lambda_{3}) = g_{H}(\lambda_{3}) = 0$,  $\lambda _3$ is a
common root of $f_H, g_H$, hence by Lemma \ref{Lemma2.1} it is
$\Re(f_{H}, g_{H}) = 0$. It is obvious that $\Re(f_{H}, g_{H})$ is a
polynomial
of $H$ with constant coefficients, therefore $H$ must be a constant.\\

\emph{Case-(ii)}: \emph{Four distinct principal curvatures}

 Now, under the assumption that shape operator
 (\ref{2.8}) has four distinct eigenvalues for $n\geq 4$, we have $\lambda_{3}=\lambda_{4}=  \ldots
 =\lambda_{n-1}$. Therefore, equation (\ref{4.1}) reduces to

\begin{equation}\label{4.19}
         (n-3)\lambda_{3}= \frac{3nH}{2} \Rightarrow \lambda_{3}=
         \frac{3nH}{2(n-3)}.
       \end{equation}

       Differentiating (\ref{4.19}) along $e_{n}$, and using (\ref{4.19}),
       T27 for $a=3$, we get

       \begin{equation}\label{4.20}
         3e_{n}(H)= H \omega_{33}^{n}.
       \end{equation}

       Again, differentiating (\ref{4.20}) along $e_{1}$, and using (\ref{4.4}), (\ref{3.23}) for $l=1$, we obtain

       \begin{equation}\label{4.21}
         H \omega_{33}^{n} \omega_{33}^{1}=0,
       \end{equation}
       which implies that either $H=0$ or $\omega_{33}^{n}=0$ or
       $\omega_{33}^{1}=0$. Now, $\omega_{33}^{n}$ can not be zero, since
       (\ref{4.20}) gives $e_{n}(H)=0$ which is a
       contradiction. Therefore, we have either $H=0$ or
       $\omega_{33}^{1}=0$.

       In similar way, if we differentiate (\ref{4.20}) along $e_{2}$, and using (\ref{4.4}), (\ref{3.23}) for $l=2$, we have either $H=0$ or
       $\omega_{33}^{2}=0$.

       Now, we claim that $H=0$. If $H$ is non-zero then we have
       $\omega_{33}^{1}=\omega_{33}^{2}=0$ and using Lemma \ref{Lemma3.4} for
       $a=3$, (\ref{4.19}), we get

       \begin{equation}\label{4.22}
         \frac{3nH}{2(n-3)}\mu=0.
       \end{equation}

       Since $\mu$ is non-zero so $H=0$, which proves our claim.
       Hence $M^{n}_{1}$ is minimal.\\

       \emph{Case-(iii)}: \emph{Three distinct principal curvatures}

       Now, for shape operator
      (\ref{2.8}) having three distinct eigenvalues for $n\geq 3$, we have $\lambda_{3}=\lambda_{4}=  \ldots
        =\lambda_{n}= \frac{-nH}{2}$. Therefore, from (\ref{3.8}), we
        $e_{n}(H)=0$ which is a contradiction.\\

         \emph{Case-(iv)}: \emph{Two distinct principal curvatures}

         If the shape operator (\ref{2.8}) has two distinct complex eigenvalues for $n=2$ then there is nothing to prove.\\

         Hence, our main Theorem \ref{Theorem1.1} follows from cases $(i)$, $(ii)$, $(iii)$ and $(iv)$.\\

\section{\textbf{Biconservative Lorentz hypersurfaces with constant length of second fundamental form}}

In this section, we study biconservative Lorentz hypersurface in
$E_{1}^{n+1}$ having shape operator (\ref{2.8}) for six distinct
eigenvalues with
constant length of second fundamental form.\\

$$\textbf{Proof of Theorem 1.2}$$
 Under the assumption that shape operator
 (\ref{2.8}) has six distinct eigenvalues, we can consider $\lambda_{3}=\lambda_{4}= \ldots =\lambda_{r}$, $\lambda_{r+1}=\lambda_{r+2}= \ldots
=\lambda_{r+s}$ and $\lambda_{r+s+1}=\lambda_{r+s+2}= \ldots
 =\lambda_{n-1}$.

 So, using (\ref{2.4}), we have trace$S=nH$ and using (\ref{3.21}) and the value
of $\lambda_{n}=-\frac{nH}{2}$, we get

 \begin{equation}\label{5.1}
 (r-2)\lambda_{\widetilde{A}}+s\lambda_{\widetilde{B}}+(n-r-s-1)\lambda_{\widetilde{C}}= \frac{3nH}{2}
 \end{equation}

 And, expressions (\ref{3.1}) reduce to
\begin{equation}\label{5.2}
\begin{array}{rcl}
S(e_{1})=\mu e_{2},\hspace{.2 cm} S(e_{2})=-\mu e_{1},\hspace{.2 cm}
S(e_{\widetilde{A}})=\lambda_{\widetilde{A}}
e_{\widetilde{A}},\\\hspace{.2 cm}
S(e_{\widetilde{B}})=\lambda_{\widetilde{B}} e_{\widetilde{B}},
\hspace{.2 cm} S(e_{\widetilde{C}})=\lambda_{\widetilde{C}}
e_{\widetilde{C}}, S(e_{n})=\lambda_{n} e_{n},
\end{array}
 \end{equation}
where $\widetilde{A}=3,4, \dots ,r$, $\widetilde{B}= r+1, r+2, \dots
,r+s$ and $\widetilde{C}= r+s+1, r+s+2, \dots ,n-1$.

Now, the hypersurface $M^{n}_{1}$ has second fundamental form $h$ of
constant length. So, we can write $||h||^{2}$=trace$S^{2}$ = $k_{1}$
where $k_{1}$ is the constant.

Therefore, we have
 \begin{equation}\label{5.3}
 (r-2)\lambda_{\widetilde{A}}^{2}+s\lambda_{\widetilde{B}}^{2}+(n-r-s-1)\lambda_{\widetilde{C}}^{2}-2\mu^{2}=k_{1}-
 \frac{n^{2}H^{2}}{4}.
 \end{equation}

Differentiating (\ref{5.1}) along $e_{n}$ and using T27 for
$a=\widetilde{A},\widetilde{B},\widetilde{C}$, we get
\begin{equation}\label{5.4}
3ne_{n}(H)=(r-2)(2\lambda_{\widetilde{A}}+nH)\omega^{n}_{\widetilde{A}\widetilde{A}}+s(2\lambda_{\widetilde{B}}+nH)\omega^{n}_{\widetilde{B}\widetilde{B}}
+(n-r-s-1)(2\lambda_{\widetilde{C}}+nH)\omega^{n}_{\widetilde{C}\widetilde{C}},
\end{equation}

Again, differentiating (\ref{5.3}) along $e_{n}$ and using T27 for
$a=\widetilde{A},\widetilde{B},\widetilde{C},$ Lemma \ref{Lemma3.3}
and (\ref{5.4}), we get
\begin{equation}\label{5.5}
\begin{array}{rcl}
(r-2)(2\lambda_{\widetilde{A}}+nH)(6\lambda_{\widetilde{A}}+nH)\omega^{n}_{\widetilde{A}\widetilde{A}}
+s(2\lambda_{\widetilde{B}}+nH)(6\lambda_{\widetilde{B}}+nH)\omega^{n}_{\widetilde{B}\widetilde{B}}\\
+(n-r-s-1)(2\lambda_{\widetilde{C}}+nH)(6\lambda_{\widetilde{C}}+nH)\omega^{n}_{\widetilde{C}\widetilde{C}}=0,
\end{array}
\end{equation}

Differentiating (\ref{5.1}) along $e_{1}$ and $e_{2}$,
alternatively, we obtain
\begin{equation}\label{5.6}
 (r-2)e_{1}(\lambda_{\widetilde{A}})+se_{1}(\lambda_{\widetilde{B}})+(n-r-s-1)e_{1}(\lambda_{\widetilde{C}})=0,
 \end{equation}

 \begin{equation}\label{5.7}
 (r-2)e_{2}(\lambda_{\widetilde{A}})+se_{2}(\lambda_{\widetilde{B}})+(n-r-s-1)e_{2}(\lambda_{\widetilde{C}})=0.
 \end{equation}

 Now, differentiating (\ref{5.3}) along $e_{1}$ and eliminating
 $e_{1}(\lambda_{\widetilde{A}})$ using (\ref{5.6}), we get
\begin{equation}\label{5.8}
 se_{1}(\lambda_{\widetilde{B}})(\lambda_{\widetilde{B}}-\lambda_{\widetilde{A}})+(n-r-s-1)e_{1}(\lambda_{\widetilde{C}})(\lambda_{\widetilde{C}}-\lambda_{\widetilde{A}})=0.
 \end{equation}

Similarly, differentiating (\ref{5.3}) along $e_{2}$ and eliminating
 $e_{2}(\lambda_{\widetilde{A}})$ using (\ref{5.7}), we get
\begin{equation}\label{5.9}
 se_{2}(\lambda_{\widetilde{B}})(\lambda_{\widetilde{B}}-\lambda_{\widetilde{A}})+(n-r-s-1)e_{2}(\lambda_{\widetilde{C}})(\lambda_{\widetilde{C}}-\lambda_{\widetilde{A}})=0.
 \end{equation}

 Now, using T7, T15 for $a=\widetilde{B},\widetilde{C}$ in (\ref{5.8}), (\ref{5.9}) and using (\ref{3.5}) we obtain
\begin{equation}\label{5.10}
 s(\lambda_{\widetilde{B}}\omega^{1}_{\widetilde{B}\widetilde{B}}-\mu \omega^{2}_{\widetilde{B}\widetilde{B}})(\lambda_{\widetilde{B}}-\lambda_{\widetilde{A}})
 +(n-r-s-1)(\lambda_{\widetilde{C}}\omega^{1}_{\widetilde{C}\widetilde{C}}-\mu
 \omega^{2}_{\widetilde{C}\widetilde{C}})(\lambda_{\widetilde{C}}-\lambda_{\widetilde{A}})=0,
 \end{equation}
 \begin{equation}\label{5.11}
 s(\lambda_{\widetilde{B}}\omega^{2}_{\widetilde{B}\widetilde{B}}+\mu \omega^{1}_{\widetilde{B}\widetilde{B}})
 (\lambda_{\widetilde{B}}-\lambda_{\widetilde{A}})+(n-r-s-1)(\lambda_{\widetilde{C}}\omega^{2}_{\widetilde{C}\widetilde{C}}+\mu
 \omega^{1}_{\widetilde{C}\widetilde{C}})(\lambda_{\widetilde{C}}-\lambda_{\widetilde{A}})=0,
 \end{equation}
 respectively.

 Now, solving (\ref{5.10}) and (\ref{5.11}) for $\omega^{1}_{\widetilde{B}\widetilde{B}}$ and
 $\omega^{2}_{\widetilde{B}\widetilde{B}}$, we find
\begin{equation}\label{5.12}
P_{1}\omega^{2}_{\widetilde{B}\widetilde{B}}=-Q_{1}
\omega^{1}_{\widetilde{C}\widetilde{C}} -
R_{1}\omega^{2}_{\widetilde{C}\widetilde{C}},
\end{equation}

\begin{equation}\label{5.13}
P_{1}\omega^{1}_{\widetilde{B}\widetilde{B}}=Q_{1}
\omega^{2}_{\widetilde{C}\widetilde{C}} -
R_{1}\omega^{1}_{\widetilde{C}\widetilde{C}},
\end{equation}
where
$P_{1}=s(\lambda_{\widetilde{B}}-\lambda_{\widetilde{A}})(\lambda_{\widetilde{B}}^{2}+\mu^{2}),$
$Q_{1}=\mu(n-r-s-1)(\lambda_{\widetilde{C}}-\lambda_{\widetilde{A}})(\lambda_{\widetilde{B}}-\lambda_{\widetilde{C}})$,
$R_{1}=(n-r-s-1)(\lambda_{\widetilde{C}}-\lambda_{\widetilde{A}})(\lambda_{\widetilde{B}}\lambda_{\widetilde{C}}+\mu^{2})$

 Now, differentiating (\ref{5.3}) along $e_{1}$ and eliminating
 $e_{1}(\lambda_{\widetilde{B}})$ using (\ref{5.6}), we get
\begin{equation}\label{5.14}
 (r-2)e_{1}(\lambda_{\widetilde{A}})(\lambda_{\widetilde{A}}-\lambda_{\widetilde{B}})+(n-r-s-1)e_{1}(\lambda_{\widetilde{C}})(\lambda_{\widetilde{C}}-\lambda_{\widetilde{B}})=0.
 \end{equation}

Similarly, differentiating (\ref{5.3}) along $e_{2}$ and eliminating
 $e_{2}(\lambda_{\widetilde{B}})$ using (\ref{5.7}), we get
\begin{equation}\label{5.15}
 (r-2)e_{2}(\lambda_{\widetilde{A}})(\lambda_{\widetilde{A}}-\lambda_{\widetilde{B}})+(n-r-s-1)e_{2}(\lambda_{\widetilde{C}})(\lambda_{\widetilde{C}}-\lambda_{\widetilde{B}})=0.
 \end{equation}

 Now, using T7, T15 for $a=\widetilde{A},\widetilde{C}$ in (\ref{5.14}), (\ref{5.15}) and using (\ref{3.5}) we obtain
\begin{equation}\label{5.16}
 (r-2)(\lambda_{A}\omega^{1}_{\widetilde{A}\widetilde{A}}-\mu \omega^{2}_{\widetilde{A}\widetilde{A}})(\lambda_{\widetilde{A}}-\lambda_{\widetilde{B}})
 +(n-r-s-1)(\lambda_{\widetilde{C}}\omega^{1}_{\widetilde{C}\widetilde{C}}-\mu
 \omega^{2}_{\widetilde{C}\widetilde{C}})(\lambda_{\widetilde{C}}-\lambda_{\widetilde{B}})=0,
 \end{equation}
 \begin{equation}\label{5.17}
 (r-2)(\lambda_{\widetilde{A}}\omega^{2}_{\widetilde{A}\widetilde{A}}+\mu \omega^{1}_{\widetilde{A}\widetilde{A}})
 (\lambda_{\widetilde{A}}-\lambda_{\widetilde{B}})+(n-r-s-1)(\lambda_{\widetilde{C}}\omega^{2}_{\widetilde{C}\widetilde{C}}+\mu
 \omega^{1}_{\widetilde{C}\widetilde{C}})(\lambda_{\widetilde{C}}-\lambda_{\widetilde{B}})=0,
 \end{equation}
 respectively.

 Now, solving (\ref{5.16}) and (\ref{5.17}) for $\omega^{1}_{\widetilde{A}\widetilde{A}}$ and
 $\omega^{2}_{\widetilde{A}\widetilde{A}}$, we find
\begin{equation}\label{5.18}
P_{2}\omega^{2}_{\widetilde{A}\widetilde{A}}=-Q_{2}
\omega^{1}_{\widetilde{C}\widetilde{C}} -
R_{2}\omega^{2}_{\widetilde{C}\widetilde{C}},
\end{equation}

\begin{equation}\label{5.19}
P_{2}\omega^{1}_{\widetilde{A}\widetilde{A}}=Q_{2}
\omega^{2}_{\widetilde{C}\widetilde{C}} -
R_{2}\omega^{1}_{\widetilde{C}\widetilde{C}},
\end{equation}
where
$P_{2}=s(\lambda_{\widetilde{B}}-\lambda_{\widetilde{A}})(\lambda_{\widetilde{A}}^{2}+\mu^{2}),$
$Q_{2}=\mu(n-r-s-1)(\lambda_{\widetilde{A}}-\lambda_{\widetilde{C}})(\lambda_{\widetilde{B}}-\lambda_{\widetilde{C}})$,
$R_{2}=(n-r-s-1)(\lambda_{\widetilde{B}}-\lambda_{\widetilde{C}})(\lambda_{\widetilde{A}}\lambda_{\widetilde{C}}+\mu^{2}).$\\

Using (\ref{5.12}), (\ref{5.13}), (\ref{5.18}) and (\ref{5.19}), we
can conclude the following:
\begin{lem}\label{Lemma5.1}
 Let $M^{n}_{1}$ be biconservative hypersurface in $E^{n+1}_{1}$ with constant length of second fundamental form having the shape operator (\ref{2.8})
 with six distinct principal curvatures with respect to suitable orthonormal basis  $\{e_{1},
e_{2}, ..., e_{n}\}$. If grad$H$ in the direction of $e_{n }$, then
we have

\begin{center}
$(\omega^{1}_{\widetilde{A}\widetilde{A}})^{2}+(\omega^{2}_{\widetilde{A}\widetilde{A}})^{2}=\frac{1}{P_{2}^{2}}(Q_{2}^{2}+R_{2}^{2})[(\omega^{1}_{\widetilde{C}\widetilde{C}})^{2}+(\omega^{2}_{CC})^{2}]$
\end{center}
\begin{center}
$(\omega^{1}_{\widetilde{B}\widetilde{B}})^{2}+(\omega^{2}_{\widetilde{B}\widetilde{B}})^{2}=\frac{1}{P_{1}^{2}}(Q_{1}^{2}+R_{1}^{2})[(\omega^{1}_{\widetilde{C}\widetilde{C}})^{2}+(\omega^{2}_{\widetilde{C}\widetilde{C}})^{2}]$
\end{center}
\begin{center}
$\omega^{1}_{\widetilde{A}\widetilde{A}}\omega^{1}_{\widetilde{C}\widetilde{C}}+\omega^{2}_{\widetilde{A}\widetilde{A}}\omega^{2}_{\widetilde{C}\widetilde{C}}=-\frac{R_{2}}{P_{2}}[(\omega^{1}_{\widetilde{C}\widetilde{C}})^{2}+(\omega^{2}_{\widetilde{C}\widetilde{C}})^{2}]$
\end{center}
\begin{center}
$\omega^{2}_{\widetilde{A}\widetilde{A}}\omega^{1}_{\widetilde{C}\widetilde{C}}-\omega^{1}_{\widetilde{A}\widetilde{A}}\omega^{2}_{\widetilde{C}\widetilde{C}}=-\frac{Q_{2}}{P_{2}}[(\omega^{1}_{\widetilde{C}\widetilde{C}})^{2}+(\omega^{2}_{\widetilde{C}\widetilde{C}})^{2}]$
\end{center}
\begin{center}
$\omega^{1}_{\widetilde{B}\widetilde{B}}\omega^{1}_{\widetilde{C}\widetilde{C}}+\omega^{2}_{\widetilde{B}\widetilde{B}}\omega^{2}_{\widetilde{C}\widetilde{C}}=-\frac{R_{1}}{P_{1}}[(\omega^{1}_{\widetilde{C}\widetilde{C}})^{2}+(\omega^{2}_{\widetilde{C}\widetilde{C}})^{2}]$
\end{center}
\begin{center}
$\omega^{2}_{\widetilde{B}\widetilde{B}}\omega^{1}_{\widetilde{C}\widetilde{C}}-\omega^{1}_{\widetilde{B}\widetilde{B}}\omega^{2}_{\widetilde{C}\widetilde{C}}=-\frac{Q_{1}}{P_{1}}[(\omega^{1}_{\widetilde{C}\widetilde{C}})^{2}+(\omega^{2}_{\widetilde{C}\widetilde{C}})^{2}]$
\end{center}
\begin{center}
$\omega^{1}_{\widetilde{B}\widetilde{B}}\omega^{1}_{\widetilde{A}\widetilde{A}}+\omega^{2}_{\widetilde{B}\widetilde{B}}\omega^{2}_{\widetilde{A}\widetilde{A}}=\frac{1}{P_{1}P_{2}}(Q_{1}Q_{2}+R_{1}R_{2})[(\omega^{1}_{\widetilde{C}\widetilde{C}})^{2}+(\omega^{2}_{\widetilde{C}\widetilde{C}})^{2}]$
\end{center}
\begin{center}
$\omega^{2}_{\widetilde{A}\widetilde{A}}\omega^{1}_{\widetilde{B}\widetilde{B}}-\omega^{1}_{\widetilde{A}\widetilde{A}}\omega^{2}_{\widetilde{B}\widetilde{B}}=\frac{1}{P_{1}P_{2}}(Q_{2}R_{1}-R_{2}Q_{1})[(\omega^{1}_{\widetilde{C}\widetilde{C}})^{2}+(\omega^{2}_{\widetilde{C}\widetilde{C}})^{2}]$
\end{center}
\end{lem}

 \vspace{.5 cm}
Now, differentiating (\ref{5.4}) along $e_{1}$, $e_{2}$
alternatively and using Lemma \ref{Lemma3.2}, T7, T15 for
$a=\widetilde{A},\widetilde{B},\widetilde{C}$, (\ref{4.4}),
(\ref{3.21}) and (\ref{3.5}), we obtain
\begin{equation}\label{5.20}
\begin{array}{rcl}
(r-2)[-2\mu
\omega^{2}_{\widetilde{A}\widetilde{A}}+(4\lambda_{\widetilde{A}}+nH)\omega^{1}_{\widetilde{A}\widetilde{A}}]\omega^{n}_{\widetilde{A}\widetilde{A}}+s[-2\mu
\omega^{2}_{\widetilde{B}\widetilde{B}}+(4\lambda_{\widetilde{B}}+nH)\omega^{1}_{\widetilde{B}\widetilde{B}}]\omega^{n}_{\widetilde{B}\widetilde{B}}\\
+(n-r-s-1)[-2\mu
\omega^{2}_{\widetilde{C}\widetilde{C}}+(4\lambda_{\widetilde{C}}+nH)\omega^{1}_{\widetilde{C}\widetilde{C}}]\omega^{n}_{\widetilde{C}\widetilde{C}}=0,
\end{array}
\end{equation}
\begin{equation}\label{5.21}
\begin{array}{rcl}
(r-2)[2\mu
\omega^{1}_{\widetilde{A}\widetilde{A}}+(4\lambda_{\widetilde{A}}+nH)\omega^{2}_{\widetilde{A}\widetilde{A}}]\omega^{n}_{\widetilde{A}\widetilde{A}}+s[2\mu
\omega^{1}_{\widetilde{B}\widetilde{B}}+(4\lambda_{\widetilde{B}}+nH)\omega^{2}_{\widetilde{B}\widetilde{B}}]\omega^{n}_{\widetilde{B}\widetilde{B}}\\
+(n-r-s-1)[2\mu
\omega^{1}_{\widetilde{C}\widetilde{C}}+(4\lambda_{\widetilde{C}}+nH)\omega^{2}_{\widetilde{C}\widetilde{C}}]\omega^{n}_{\widetilde{C}\widetilde{C}}=0,
\end{array}
\end{equation}
respectively.

Now, (\ref{5.5}), (\ref{5.20}), (\ref{5.21}) form a homogeneous
system of equations in $\omega^{n}_{\widetilde{A}\widetilde{A}}$,
$\omega^{n}_{\widetilde{B}\widetilde{B}}$ and
$\omega^{n}_{\widetilde{C}\widetilde{C}}$ having non trivial
solution. Therefore, discriminant $D=0$ which gives
\begin{equation}\label{5.22}
\begin{array}{lcl}
(2\lambda_{\widetilde{A}}+nH)(6\lambda_{\widetilde{A}}+nH)\{8\mu(\lambda_{\widetilde{B}}-\lambda_{\widetilde{C}})
(\omega^{1}_{\widetilde{B}\widetilde{B}}\omega^{1}_{\widetilde{C}\widetilde{C}}+\omega^{2}_{\widetilde{B}\widetilde{B}}\omega^{2}_{\widetilde{C}\widetilde{C}})
\\+[(4\lambda_{\widetilde{C}}+nH)(4\lambda_{\widetilde{B}}+nH)+4\mu^{2}](\omega^{1}_{\widetilde{B}\widetilde{B}}\omega^{2}_{\widetilde{C}\widetilde{C}}-\omega^{2}_{\widetilde{B}\widetilde{B}}\omega^{1}_{\widetilde{C}\widetilde{C}})\}\\-
(2\lambda_{\widetilde{B}}+nH)(6\lambda_{\widetilde{B}}+nH)\{8\mu(\lambda_{\widetilde{A}}-\lambda_{\widetilde{C}})(\omega^{1}_{\widetilde{A}\widetilde{A}}\omega^{1}_{\widetilde{C}\widetilde{C}}+\omega^{2}_{\widetilde{A}\widetilde{A}}\omega^{2}_{\widetilde{C}\widetilde{C}})
\\+[(4\lambda_{\widetilde{C}}+nH)(4\lambda_{\widetilde{A}}+nH)+4\mu^{2}](\omega^{1}_{\widetilde{A}\widetilde{A}}\omega^{2}_{\widetilde{C}\widetilde{C}}
-\omega^{2}_{\widetilde{A}\widetilde{A}}\omega^{1}_{\widetilde{C}\widetilde{C}})\}\\+
(2\lambda_{\widetilde{C}}+nH)(6\lambda_{\widetilde{C}}+nH)\{8\mu(\lambda_{\widetilde{A}}-\lambda_{\widetilde{B}})(\omega^{1}_{\widetilde{A}\widetilde{A}}\omega^{1}_{\widetilde{B}\widetilde{B}}+\omega^{2}_{\widetilde{A}\widetilde{A}}\omega^{2}_{\widetilde{B}\widetilde{B}})
\\+[(4\lambda_{\widetilde{A}}+nH)(4\lambda_{\widetilde{B}}+nH)+4\mu^{2}](\omega^{1}_{\widetilde{A}\widetilde{A}}\omega^{2}_{\widetilde{B}\widetilde{B}}
-\omega^{2}_{\widetilde{A}\widetilde{A}}\omega^{1}_{\widetilde{B}\widetilde{B}})\}=0,
\end{array}
\end{equation}
which by using Lemma \ref{Lemma5.1} and the values of $P_{1}$,
$P_{2}$, $Q_{1}$, $Q_{2}$, $R_{1}$, $R_{2}$ reduces to
\begin{equation}\label{5.23}
\mu(\lambda_{\widetilde{A}}-\lambda_{\widetilde{B}})(\lambda_{\widetilde{B}}-\lambda_{\widetilde{C}})(\lambda_{\widetilde{C}}-\lambda_{\widetilde{A}})[(\omega_{\widetilde{C}\widetilde{C}}^{1})^{2}+
(\omega_{\widetilde{C}\widetilde{C}}^{2})^{2}]f(\lambda_{\widetilde{A}},\lambda_{\widetilde{B}},\lambda_{\widetilde{C}},H)=0
\end{equation}
where

 $\begin{array}{rcl}
f(\lambda_{\widetilde{A}},\lambda_{\widetilde{B}},\lambda_{\widetilde{C}},H)=
(r-2)(2\lambda_{\widetilde{A}}+nH)(6\lambda_{\widetilde{A}}+nH)(\lambda_{\widetilde{A}}^{2}+\mu^{2})\{4\mu^{2}-8\lambda_{\widetilde{B}}\lambda_{\widetilde{C}}
-4nH(\lambda_{\widetilde{B}}\\+\lambda_{\widetilde{C}})
-n^{2}H^{2}\}+s(2\lambda_{\widetilde{B}}+nH)(6\lambda_{\widetilde{B}}+nH)(\lambda_{\widetilde{C}}^{2}+\mu^{2})
\{4\mu^{2}-8\lambda_{\widetilde{A}}\lambda_{\widetilde{C}}-4nH(\lambda_{\widetilde{A}}+\lambda_{\widetilde{C}})
-\\n^{2}H^{2}\}+(n-r-s-1)(2\lambda_{\widetilde{C}}+nH)(6\lambda_{\widetilde{C}}+nH)(\lambda_{\widetilde{B}}^{2}
+\mu^{2})
\{4\mu^{2}-8\lambda_{\widetilde{B}}\lambda_{\widetilde{A}}-4nH\\(\lambda_{\widetilde{B}}+
\lambda_{\widetilde{A}}) -n^{2}H^{2}\}.
\end{array}$

\vspace{.3 cm}

 Now, we claim that
$(\omega_{\widetilde{C}\widetilde{C}}^{1})^{2}+(\omega_{\widetilde{C}\widetilde{C}}^{2})^{2}\neq
       0$.\\

      If $(\omega_{\widetilde{C}\widetilde{C}}^{1})^{2}+(\omega_{\widetilde{C}\widetilde{C}}^{2})^{2}=0$, we have,
      $\omega_{\widetilde{C}\widetilde{C}}^{1}=\omega_{\widetilde{C}\widetilde{C}}^{2}=0$ as connection coefficients are real numbers. Then, using (\ref{5.12}),
       (\ref{5.13}), (\ref{5.18}) and (\ref{5.19}),
      we have $\omega_{\widetilde{A}\widetilde{A}}^{1}=\omega_{\widetilde{A}\widetilde{A}}^{2}=\omega_{\widetilde{B}\widetilde{B}}^{1}=\omega_{\widetilde{B}\widetilde{B}}^{2}=0$.

      Therefore, using Lemma \ref{Lemma3.4} for $a=\widetilde{A}, \widetilde{B}, \widetilde{C},$ we obtain

      \begin{equation}\label{5.24}
         \lambda_{\widetilde{A}}\mu=0, \quad \lambda_{\widetilde{B}}\mu=0, \quad \lambda_{\widetilde{C}}\mu=0
       \end{equation}
       which implies $\lambda_{\widetilde{A}}=\lambda_{\widetilde{B}}=\lambda_{\widetilde{C}}=0$.
       Now, using T27 for $a = \widetilde{A}, \widetilde{B}, \widetilde{C}$, we obtain
       $\omega_{\widetilde{A}\widetilde{A}}^{n}=\omega_{\widetilde{B}\widetilde{B}}^{n}=\omega_{\widetilde{C}\widetilde{C}}^{n}=0$ and (\ref{5.4}) gives $e_{n}(H)=0$
        which is a contradiction. Hence our claim is proved.

     Now, using the fact that $\lambda_{\widetilde{A}},\lambda_{\widetilde{B}},\lambda_{\widetilde{C}}$
     are distinct and $\mu$ is non zero, (\ref{5.23}) implies that
     \begin{equation}\label{5.25}
         f(\lambda_{\widetilde{A}}, \lambda_{\widetilde{B}}, \lambda_{\widetilde{C}}, H)=0
       \end{equation}

       Now, eliminating $\lambda_{\widetilde{A}}$ from (\ref{5.25}) using (\ref{5.1}), we
       get
       \begin{equation}\label{5.26}
       \begin{array}{lcl}
       f_{1}(\lambda_{\widetilde{B}}, \lambda_{\widetilde{C}}, H)=[nH(r+1)-2s\lambda_{\widetilde{B}}-2(n-r-s-1)\lambda_{\widetilde{C}}][nH(r+7)-6s\lambda_{\widetilde{B}}\\ -6(n-r-s-1)\lambda_{\widetilde{C}}]
       [9n^{2}H^{2}+4s^{2}\lambda_{\widetilde{B}}^{2}+4(n-r-s-1)^{2}\lambda_{\widetilde{C}}^{2}\\ -12nH(s\lambda_{\widetilde{B}}+(n-r-s-1)\lambda_{\widetilde{C}})-8s(n-r-s-1)\lambda_{\widetilde{B}}\lambda_{\widetilde{C}}\\
       +4\mu^{2}(r-2)^{2}][4\mu^{2}-8\lambda_{\widetilde{A}}\lambda_{\widetilde{C}}-4nH(\lambda_{\widetilde{A}}+\lambda_{\widetilde{C}})
       -n^{2}H^{2}]+\\ 4s(r-2)^{2}(2\lambda_{\widetilde{B}}+nH)(6\lambda_{\widetilde{B}}+nH)(\lambda_{\widetilde{C}}^{2}+\mu^{2})[4\mu^{2}(r-2)-\\
       4\lambda_{\widetilde{C}}(3nH-2s\lambda_{\widetilde{B}}-2(n-r-s-1)\lambda_{\widetilde{C}})-2nH\{3nH-2s\lambda_{\widetilde{B}}\\ -2(n-2r-s+1)\lambda_{\widetilde{C}}\}-(r-2)n^{2}H^{2}]+
       4(n-r-s-1)(r-2)^{2}(2\lambda_{\widetilde{C}}\\ +nH)(6\lambda_{\widetilde{C}}+nH)(\lambda_{\widetilde{C}}^{2}+\mu^{2})[4\mu^{2}(r-2)-
       4\lambda_{\widetilde{B}}(3nH-2s\lambda_{\widetilde{B}}-\\ 2(n-r-s-1)\lambda_{\widetilde{C}})-2nH\{3nH-2(s-r+2)\lambda_{\widetilde{B}}-2(n-r-s-1)\lambda_{\widetilde{C}}\}\\
       -(r-2)n^{2}H^{2}]=0
       \end{array}
       \end{equation}

       Similarly, eliminating $\lambda_{\widetilde{A}}$ from (\ref{5.3}) using (\ref{5.1}), we
       obtain
        \begin{equation}\label{5.27}
        \begin{array}{lcl}
        g_{1}(\lambda_{\widetilde{B}}, \lambda_{\widetilde{C}}, H)=[3nH-2s\lambda_{\widetilde{B}}-2(n-r-s-1)\lambda_{\widetilde{C}}]^{2}+4s(r-2)\lambda_{\widetilde{B}}^{2}\\+4(n-r-s-1)(r-2)\lambda_{\widetilde{C}}^{2}
        -4(r-2)2\mu^{2}-4(r-2)k_{1}+ 4(r-2)\frac{n^{2}H^{2}}{4}=0
        \end{array}
        \end{equation}

        We rewrite $f_{1}(\lambda_{\widetilde{B}},\lambda_{\widetilde{C}},H)$, $g_{1}(\lambda_{\widetilde{B}},\lambda_{\widetilde{C}},H)$ as polynomials
$f_{1(H,\lambda_{\widetilde{C}})}(\lambda_{\widetilde{B}}),
g_{1(H,\lambda_{\widetilde{C}})}(\lambda_{\widetilde{B}})$ of
$\lambda_{\widetilde{B}}$ with coefficients in polynomial ring
$R_{1}[\lambda_{\widetilde{C}},H]$ over real field $\mathbb{R}$.
According to the Lemma \ref{Lemma2.1}, the equations
$f_{1(H,\lambda_{\widetilde{C}})}(\lambda_{\widetilde{B}}) = 0$ and
$g_{1(H,\lambda_{\widetilde{C}})}(\lambda_{\widetilde{B}}) = 0$ have
a common root if and only if resultant
$\Re(f_{1(H,\lambda_{\widetilde{C}})},
g_{1(H,\lambda_{\widetilde{C}})}) = 0$. It is obvious that
$\Re(f_{1(H,\lambda_{\widetilde{C}})},
g_{1(H,\lambda_{\widetilde{C}})})$ is a polynomial of
$\lambda_{\widetilde{C}}$ and $H$. So, we have
\begin{equation}\label{5.28}
f_{2}(\lambda_{\widetilde{C}},H)=\Re(f_{1(H,\lambda_{\widetilde{C}})},
g_{1(H,\lambda_{\widetilde{C}})}) = 0
\end{equation}

       Differentiating (\ref{5.28}) along $e_{1}$ and $e_{2}$, we have

       \begin{equation}\label{5.29}
         e_{1}(\lambda_{\widetilde{C}})(g_{2}(\lambda_{\widetilde{C}},H))=0,
       \end{equation}
       and
\begin{equation}\label{5.30}
         e_{2}(\lambda_{\widetilde{C}})(g_{2}(\lambda_{\widetilde{C}},H))=0,
       \end{equation}
       respectively.

Now, if $g_{2}(\lambda_{\widetilde{C}},H)\neq 0$, we have
$e_{1}(\lambda_{\widetilde{C}})=0$ and
$e_{2}(\lambda_{\widetilde{C}})=0$ which implies from T7, T15 for
$a=\widetilde{C}$, (\ref{5.12}), (\ref{5.13}), (\ref{5.18}) and
(\ref{5.19}) that
$\omega_{\widetilde{C}\widetilde{C}}^{1}=\omega_{\widetilde{C}\widetilde{C}}^{2}=\omega_{\widetilde{B}\widetilde{B}}^{1}=\omega_{\widetilde{B}\widetilde{B}}^{2}
=\omega_{\widetilde{A}\widetilde{A}}^{1}=\omega_{\widetilde{A}\widetilde{A}}^{2}=0$.
And, we have already proved from (\ref{5.24}) that it arises to a
contradiction.

Therefore, we have

\begin{equation}\label{5.31}
         g_{2}(\lambda_{\widetilde{C}},H)= 0,
       \end{equation}
  which is also a polynomial equation in $\lambda_{\widetilde{C}}$ and $H$ of degree $2$.

Again, we rewrite $f_{2}(\lambda_{\widetilde{C}},H)$,
$g_{2}(\lambda_{\widetilde{C}},H)$ as polynomials
$f_{2(H)}(\lambda_{\widetilde{C}}),
g_{2(H)}(\lambda_{\widetilde{C}})$ of $\lambda_{\widetilde{C}}$ with
coefficients in polynomial ring $R_{2}[H]$ over real field
$\mathbb{R}$. According to the Lemma \ref{Lemma2.1}, the equations
$f_{2(H)}(\lambda_{\widetilde{C}}) = 0$ and
$g_{2(H)}(\lambda_{\widetilde{C}}) = 0$ have a common root if and
only if $\Re(f_{2(H)}, g_{2(H)}) = 0$. It is obvious that
$\Re(f_{2(H)}, g_{2(H)})$ is a polynomial of $H$ with constant
coefficients. So, $\Re(f_{2(H)}, g_{2(H)}) = 0$ which implies that
$H$ must be a constant and hence proved Theorem \ref{Theorem1.2}. \\

\textbf{\emph{Acknowledgement}}: \small \emph{The author is grateful
to Guru Gobind Singh Indraprastha University for providing IPR
Fellowship to pursue research. The author acknowledges useful
discussions and suggestions with Dr. Ram Shankar
Gupta and Dr. Andreas Arvanitoyeorgos.}\\


\bibliography{xbib}


Author's address:\\
\\
\textbf{Deepika}\\
University School of Basic and Applied Sciences,\\
Guru Gobind Singh Indraprastha University,\\
Sector-16C, Dwarka, New Delhi-110078, India.\\
\textbf{Email:} sdeep2007@gmail.com\\
\\

\end{document}